\documentclass[12pt]{article}
\usepackage{latexsym}

\newcommand{\Aut}{\hbox{{\rm Aut}}\, }

\newcommand{\Mon}{\hbox{{\rm Mon}}\, }
\newcommand{\core}[2]{#1_{_{#2}}}
\newcommand{\nclose}[2]{#1^{^{#2}}}

\def\S{$\breve{{\rm S}}$}

\newcommand{\sz}{\scriptsize}
\newcommand{\szc}{\small}

\long\def\footinfo#1{\begingroup\def\thefootnote{\fnsymbol{footnote}}\footnote[0]{\hskip-18pt #1}\endgroup}
\def\subjclass#1{2000 {\em Mathematics Subject Classification}. #1}
\def\keywords#1{ {\em Key words and phrases}. #1}

\begin{document}

\newcounter{teocounter}
\addtocounter{teocounter}{1}

\begin{center}
{\bf\large CHIRALITY GROUPS OF MAPS AND HYPERMAPS \\}
\end{center}

\bigskip\noindent
\begin{center}
\begin{tabular}{cc}
{\szc Antonio Breda D'Azevedo}\footnotemark &
  {\szc Gareth Jones} \\[-4pt]
                     &                      \\[-8pt]
{\sz Departamento de Matematica} & {\sz Mathematics Department}\\[-4pt]
{\sz Universidade de Aveiro} & {\sz University of Southampton }\\[-4pt]
{\sz Aveiro} & {\sz Southampton SO17 1BJ}\\[-4pt]
{\sz Portugal} & {\sz United Kingdom}\\[-4pt]
&                                    \\[-4pt]

\end{tabular}

\begin{tabular}{cc}
{\szc Roman Nedela}\footnotemark & {\szc Martin \S koviera}\footnotemark \\[-4pt]
                     &                      \\[-8pt]
{\sz Katedra Matematiky} &
{\sz Department of Informatics}\\[-4pt]
 {\sz Fakulta Financi\' e} &
{\sz Comenius University}\\[-4pt]
{\sz 975 49 Bansk\' a Bystrica} & {\sz
Bratislava}\\[-4pt]
{\sz Slovakia} & {\sz Slovakia}
\end{tabular}

\addtocounter{footnote}{-2}\footnotetext{Supported in part by UI\&D {\sl Matem{\'a}tica e aplica\c{c}{\~o}es} of University of
Aveiro, through Program POCTI of FCT co-financed by the European Community fund FEDER}
\addtocounter{footnote}{1}\footnotetext{ ${}^{,3}$Supported in part by the Ministry for
Education of the Slovak Republic, grant no. $1/6132/99$}
\end{center}

\footinfo{\subjclass{05C10, 05C25.}}
\footinfo{\keywords{map, hypermap, chiral, asymmetric, chirality group, chirality index.}}

\vskip20pt

{\small {\bf Abstract.}
Although the phenomenon of chirality appears in many investigations of maps and
hypermaps no detailed study of chirality seems to have been carried out.
Chirality of maps and hypermaps is not merely a binary invariant but can be
quantified by two new invariants --- the chirality group and the chirality index, the latter being
the size of the chirality group. A detailed investigation of the chirality
groups of maps and hypermaps will be the main objective of this paper.
The most extreme type of chirality arises when the chirality group coincides with the monodromy group.
Such hypermaps are  called totally chiral. Examples
of them are constructed by considering appropriate ``asymmetric'' pairs of generators for some
non-abelian simple groups. We also show that every finite abelian group is the chirality group
of some hypermap, whereas many non-abelian groups, including
symmetric and dihedral groups, cannot arise as chirality groups.}

\vskip10pt

\section{Introduction}

The word chirality (meaning handedness) was introduced by William Thomson,
better known as Lord
Kelvin \cite{KEL}, as a property of a geometrical figure which occurs when
``its image
in a plane mirror, ideally realized, cannot be brought into coincidence
with itself''. There are
many scientific examples of this phenomenon. For instance, chemists and
biologists have discovered
molecules which exist in two distinguishable enantiomers, which are mirror
images of each
other; in some cases, one form is beneficial whereas the other is
poisonous. The phenomenon of chirality
also appears in quantum mechanics, and more generally in theoretical physics.

These facts (see also \cite{BAR,Bro,FUJ,JAN,MEZ})
provide a strong motivation for the study of chirality
in mathematical models of physical structures.
Many of these models are graphs embedded in the
plane or 3-dimensional Euclidean space, or
in more complex spaces endowed with a non-Euclidean geometry.
For instance, fullerenes (recently-discovered positively-curved carbon
structures) are modelled by
3-valent polyhedra with pentagonal and hexagonal faces, the pentagons being
separated by
hexagons; negatively-curved schwarzites arise if the pentagons are replaced
with heptagons. These
are special examples of maps.

Let us be more precise. A {\em map\/} on a surface is a cellular
decomposition of a closed
surface into $0$-cells called {\em vertices\/}, $1$-cells called {\em edges\/}
and $2$-cells called {\em faces\/}. The vertices and edges of a map form
its {\em underlying graph\/}. A map is said to be {\em orientable\/} if the
supporting surface is orientable, and is {\em oriented\/} if one of two
possible orientations of the surface has been specified.
An automorphism of a map is an
automorphism of the underlying graph which extends to a
self-homeomorphism of the underlying surface.
Map automorphisms split naturally into two classes,
orientation-preserving and orientation-reversing automorphisms.
We say that a map is {\em chiral\/} if it admits no orientation-reversing
automorphism. These
concepts extend naturally to hypermaps, generalisations of maps which we
shall explain later.

Although the phenomenon of chirality appears in many investigations of maps and
hypermaps (see Coxeter and Moser \cite{CM} for instance),
no detailed study of chirality seems to have been carried out. A
starting-point which
led us to such a project was the observation that chirality
of maps and hypermaps is not merely a binary invariant but can be
quantified. Even more surprising was
the fact that different approaches to measuring chirality lead to
equivalent definitions, giving rise
to related new invariants associated with any (hyper)map ---
the chirality group and the chirality index, the latter being
the size of the chirality group. A detailed investigation of the chirality
groups
of maps and hypermaps will be the main objective of this paper.

In the investigation of maps and hypermaps it is often convenient
to replace the topological objects with their combinatorial counterparts.
Indeed, it is well-known that a map $\cal H$ on an orientable surface can
be described
by two permutations $R$ and $L$ acting on the set of darts (directed edges,
or arcs).
The permutation $R$ cyclically permutes the darts based at the same vertex,
consistently with a
chosen orientation of the surface, while $L$ interchanges the two
oppositely directed darts
sharing the same edge. Thus the orientation of the map can be encoded by
the choice of $R$
(the other possibility being $R^{-1}$), and our maps are by definition
oriented (or ``polarized'' using the language of physicists).
By connectivity, the action of the group $\Mon({\cal H})=\langle
R,L\rangle$ is transitive on the set of
darts of $\cal H$. A hypermap is obtained if we simply relax the
requirement $L^2=1$ by allowing
$L$ to be of any order, so maps are particular examples of hypermaps.
With the above notation $\cal H$ is said to be mirror symmetric or
reflexible if the assignment
$R\mapsto R^{-1}$, $L\mapsto L^{-1}$ extends to a group automorphism of
$\langle R,L\rangle$; if it does
not extend, we say that $\cal H$ is chiral.

Clearly, any hypermap $\cal H$ covers some reflexible hypermap $\cal K$,
because in the worst case
we can take $\cal K$ to be the trivial map consisting of one dart and a
vertex attached to it.
For simplicity, let us assume that the actions of $\Mon({\cal H})$
and $\Mon({\cal K})$ are both regular. Then the covering ${\cal H}\to {\cal K}$
arises from factoring out a certain normal subgroup of $\Mon({\cal H})$.
The minimal subgroup $X({\cal H})\unlhd \Mon({\cal H})$ such that ${\cal
H}/X({\cal H})$ is
a reflexible hypermap is called the {\em chirality group} of $\cal H$. It
is straightforward
that $\cal H$ is chiral if and only if $X({\cal H})$ is nontrivial. There
is a dual approach
to the definition of the chirality group, obtained by
considering the smallest reflexible hypermap which covers $\cal H$.
It is proved in Section 3 that these two approaches are equivalent.

In general, $|X({\cal H})|\le |\Mon({\cal H})|$.
The most extreme type of chirality arises when $X({\cal H})=\Mon({\cal
H})$; such hypermaps,
called totally chiral, are studied in Section~5 where we construct examples
of them by
considering an appropriate ``asymmetric'' pair of generators for some
non-abelian simple group.
In Section 6 we show that every finite abelian group is the chirality group
of some hypermap,
whereas it is shown in Section 7 that many non-abelian groups, including
symmetric and dihedral
groups, cannot arise as chirality groups. The general problem of
characterising chirality groups
remains open.

\vskip8pt

\section{Hypermaps}

By an {\em oriented hypermap\/} we mean a triple ${\cal H}=(D,R,L)$
where $D$ is a set of darts and $R$ and $L$ are two permutations
generating a permutation group $\Mon({\cal H})=\langle
R,L\rangle$, called the {\em monodromy group} of ${\cal H}$,
acting transitively on $D$. The permutations $R$ and $L$ will be called
the {\em canonical generators\/} of $\Mon({\cal H})$.
The orbits of $R$, $L$ and $RL$ on $D$ will be called the {\em
hypervertices\/}, {\em hyperedges\/}, and {\em hyperfaces\/},
respectively.

In general, the set of darts of a hypermap may be infinite,
however, our main interest lies in {\em finite\/} hypermaps, those
where the set $D$ (equivalently, the group
$\Mon({\cal H})$) is finite. {\em All hypermaps in this paper
will be finite unless the immediate context implies otherwise.}

If $m,n$ and $k$ are the orders of the permutations $RL$, $R$,
and $L$, respectively, then the {\em type\/} of $\cal H$ the triple
$(m,n,k)$, and its {\em
Euler characteristic\/} is the number $\chi({\cal H})=|D|(1/m +1/n +1/k -1)$.

A hypermap $\cal H$ is an {\em oriented map\/} if $L$ is an involution.
Since the maps and hypermaps considered in this
paper will be all oriented, the adjective oriented will be often omitted.

Maps are an algebraic abstraction
of topological maps, that is, cellular decompositions of closed
surfaces. Hypermaps generalize the notion of maps in a natural way, but
their geometric
representations are somewhat less natural; they are described,
for example, in \cite{Cori,Wa}.

An {\em automorphism\/} of a hypermap ${\cal H}=(D,R,L)$ is a
permutation $\psi$ of $D$ commuting with both $R$ and $L$. It is
straightforward to see that the automorphism group  $\Aut(\cal H)$
of $\cal H$ acts semi-regularly on the set of darts, so that
$|\Aut(\cal H)|$ divides $|D|$. If equality
holds, the action is regular, and consequently
${\cal H}$ is called an {\em orientably regular hypermap}.

If ${\cal H}=(D,R,L)$ and ${\cal H'}=(D',R',L')$ are
hypermaps then a {\em covering\/} $\psi:{\cal H}\to {\cal
H'}$ is a mapping $\psi:D\to D'$ satisfying $\psi R=R'\psi$ and
$\psi L=L'\psi$; note that a covering is necessarily surjective. We then
say that ${\cal H}$
{\em covers\/} ${\cal H'}$, and it follows that the assignment
$R\mapsto R'$ and $L\mapsto L'$ extends to a {\em canonical epimorphism\/}
$\Mon({\cal H})\to \Mon({\cal H'})$ of the monodromy groups. If $\psi$ is
an injective
covering, we have  an {\em isomorphism\/} ${\cal H}\cong {\cal H'}$.
A covering ${\cal H}\to {\cal H'}$ is {\em smooth\/} if both hypermaps have
the same type.

Let $\Delta$ denote the free product
\[\Delta =\langle r_0,r_1,r_2\mid
{r_0}^2={r_1}^2={r_2}^2=1\rangle\] and let $ \Delta^+=\langle
r_1r_2, r_2r_0\rangle$ be its even
word subgroup. The canonical generators of $\Delta^+$ will be
denoted by $\rho=r_1r_2$ and $\lambda=r_2r_0$. Observe that the triple
${\cal U} =( \Delta^+, \rho, \lambda)$, with $\rho$ and $\lambda$ acting on
$\Delta^+$ by
the left translation, is a hypermap (clearly, an infinite one)
which we will call the {\em universal hypermap}.

For any hypermap ${\cal H} =(D,R,L)$, finite or infinite, there is an
epimorphism
$\mu: \Delta^+\to \Mon({\cal H})$ sending $\rho$ to $R$ and $\lambda$ to
$L$. Consequently, $\cal H$ can be identified with the hypermap
$(\Delta^+/H, \bar\rho, \bar\lambda\,)$ whose darts are the left cosets of the
subgroup $H \leq \Delta^+$, $H$ being the
preimage of the stabiliser of a dart in $\cal H$ under the
action of $\Mon (\cal H)$, and $\bar\rho(xH) = \rho xH$ and $\bar\lambda (xH)
= lxH$. With some abuse of notation, $\mu: {\cal U} \to (\Delta^+/H,
\bar\rho, \bar\lambda\,)\cong
\cal H$ is a hypermap covering. Thus the monodromy group of any oriented
hypermap
is a quotient of $\Delta^+$, and oriented
hypermaps correspond to subgroups of $\Delta^+$. Any subgroup
$H\leq \Delta^+$ for which $(\Delta^+/H, \bar\rho, \bar\lambda)\cong \cal H$
will be called a {\em hypermap subgroup\/} for $\cal H$.

In the following statement we summarize some well-known facts on
representations of hypermaps by hypermap subgroups (see \cite{CS}).

\vskip8pt
\noindent {\bf Lemma \theteocounter\,\,} {\sl
Let $\cal H$,
${\cal H}_1$ and ${\cal H}_2$ be oriented hypermaps (not necessarily
finite).
Then the following statements hold:
\begin{description}
\item[(i)] ${\cal H}_1$ covers ${\cal H}_2$ if and only if
there are subgroups $H_1\le H_2\le \Delta^+$
such that ${\cal H}_1\cong (\Delta^+/H_1,R_1,L_1)$ and ${\cal
H}_2\cong (\Delta^+/H_2,R_2,L_2)$ where $R_i(xH_i)=\rho xH_i$
and $L_i(xH_i)=\lambda xH_i$ for $i=1,2$;
\item[(ii)] ${\cal H}_1\cong {\cal H}_2$ if and only if the corresponding
hypermap
subgroups are conjugate in $\Delta^+$;
\item[(iii)] $\cal H$ is orientably regular if and only if there
exists a normal subgroup $N\unlhd \Delta^+$ such that ${\cal
H}\cong (\Delta^+/N, \rho N, \lambda N)$.
Thus orientably regular hypermaps correspond to normal subgroups of $\Delta^+$.
\end{description}
}
\newcounter{general}
\setcounter{general}{\theteocounter}
\addtocounter{teocounter}{1}

While automorphisms of oriented hypermaps give rise to
orientation-preserving self-homeomorphisms of the underlying
surface, a hypermap may admit external symmetries
coming from self-homeomorphisms which change the orientation
of the surface. Such symmetries are called
mirror symmetries (or reflections, or
inversions). More precisely, a permutation $\psi$ of $D$ will be
called a {\em mirror symmetry\/} of an oriented hypermap ${\cal
H}=(D,R,L)$ if $\psi R=R^{-1}\psi$ and $\psi L=L^{-1}\psi$.
An orientably regular hypermap admitting mirror symmetries is
said to be {\em reflexible\/} or {\em regular\/}, since the automorphism
group acts
regularly on flags (triples formed by a mutually
incident hypervertex, hyperedge and hyperface).
On the other hand, an orientably regular hypermap
with no mirror symmetries
will be called a {\em chiral\/} hypermap.
Observe that conjugation by $r_2$ induces an automorphism
of $\Delta^+$ inverting its generators $\rho$ and
$\lambda$. Since a hypermap subgroup $H$ of an orientably
regular hypermap
${\cal H}$ is normal in $\Delta^+$, its conjugates in $\Delta$
are $H$ and $H^{r_0}=H^{r_1}=H^{r_2}$. Let $H^r$ denote this
common conjugate $H^{r_i}$. It is straightforward to see that
$H^r$ is a hypermap subgroup of the {\em mirror image} ${\cal
H}^r=(D,R^{-1},L^{-1})$.

\section{The Chirality Group and Chirality Index of a Hypermap}

The aim of this section is to introduce invariants which in some sense
measure the mirror asymmetry of a hypermap.

Let $\cal H$ be an orientably regular hypermap with hypermap subgroup $H$,
a normal subgroup of
$\Delta^+$. Then the largest normal subgroup of $\Delta$ contained in $H$
is the
group $\core{H}{\Delta}=H\cap H^r$, and the smallest normal subgroup of
$\Delta$ containing $H$ is
the group $\nclose{H}{\Delta}=HH^r$. The corresponding hypermaps
$\core{{\cal H}}{\Delta}$ and $\nclose{{\cal H}}{\Delta}$ are respectively
the smallest reflexible
hypermap that covers $\cal H$, and the largest reflexible hypermap that is
covered by $\cal H$.
In particular, if $\cal H$ is finite, so is $\core{{\cal H}}{\Delta}$,
since the intersection of
two subgroups of finite index also has finite index.

\vskip7pt \noindent {\bf Proposition \theteocounter\,\,} {\sl The
four groups $H^{\Delta}/H,\, H/H_{\Delta},\, H^{\Delta}/H^r$ and
$H^r/H_{\Delta}$ are all isomorphic to each other.}
\vskip5pt
\newcounter{keki}
\setcounter{keki}{\theteocounter} \addtocounter{teocounter}{1}
\noindent {\sl Proof.}
The third isomorphism theorem gives
$$H^{\Delta}/H=HH^r/H\cong H^r/(H\cap H^r)=H^r/H_{\Delta},$$
and similarly $H^{\Delta}/H^r\cong H/H_{\Delta}$. Conjugation by
a generator $r_i$ of $\Delta$ induces isomorphisms
$H^{\Delta}/H\cong H^{\Delta}/H^r$ and $H/H_{\Delta}\cong
H^r/H_{\Delta}$.
\hspace{\fill} $\Box$ \vskip5pt

We will call this common group the {\em chirality group\/}
$X({\cal H})$ of $\cal H$, and its order the {\em
chirality index\/} $\kappa=\kappa({\cal H})$ of $\cal H$. Thus
$\cal H$ is reflexible if and only if $\kappa=1$, and in general $X({\cal H})$
and $\kappa({\cal H})$ can be regarded as algebraic and numerical measures
of how far $\cal H$ deviates from being reflexible.

\vskip7pt
\noindent {\bf Theorem \theteocounter\,\,} {\em Let $\cal H$ be an orientably
regular hypermap with chirality index $\kappa$. Then
${\cal H}_{\Delta}\to {\cal H}$
and ${\cal H}\to {\cal H}^{\Delta}$ are both $\kappa$-sheeted
regular coverings with covering transformation group isomorphic
to the chirality group $X({\cal H})$. Moreover, the covering ${\cal
H}_{\Delta}\to {\cal H}$ is
smooth.}
\vskip5pt
\newcounter{skfc}
\setcounter{skfc}{\theteocounter}
\addtocounter{teocounter}{1}
\noindent {\sl Proof.}
\noindent
Let $H$ be the hypermap subgroup for $\cal H$.
Since $\core{H}{\Delta}$ is normal in $H$ and $H$ is normal in
$\nclose{H}{\Delta}$, these coverings are both regular. By Proposition
\thekeki\,
the corresponding covering transformation groups
$H/H_{\Delta}$ and $H^{\Delta}/H$ are isomorphic, and they coincide
with $X({\cal H})$. The number of sheets of the covering is the index
$|H:H_{\Delta}|$ or
$|H^{\Delta}:H|$, equal in each case to $|X({\cal H})|=\kappa$.

Since $r_2^{-1}\rho r_2=\rho^{-1}$, $r_2^{-1}\lambda
r_2=\lambda^{-1}$, and
$r_2^{-1}(\rho\lambda)r_2=\lambda(\rho\lambda)^{-1}\lambda^{-1}$,
the same powers of $\rho, \lambda$ or $\rho\lambda$ lie in $H$ as
in $H^r$, and hence as in $H\cap H^r=\core{H}{\Delta}$. Thus
$\core{{\cal H}}{\Delta}$ has the same type as $\cal H$, and is
therefore a smooth covering of $\cal H$. $\Box$

\vskip7pt

In general, the covering ${\cal H}\to {\cal H}^{\Delta}$ need not be
smooth; indeed, we shall see
in Section~5 that there are nontrivial hypermaps $\cal H$ for which ${\cal
H}^{\Delta}$ is the
trivial hypermap, of type $(1,1,1)$.

\vskip7pt
\noindent {\bf Corollary \theteocounter\,\,} {\em If $\cal H$ is
an orientably regular hypermap with chirality index $\kappa$, then
$\chi(\core{{\cal H}}{\Delta})=\kappa\chi({\cal H})$.}
\newcounter{Cskfc}
\setcounter{Cskfc}{\theteocounter}
\addtocounter{teocounter}{1}

\vskip7pt
\noindent {\bf Proposition \theteocounter\,\,}
{\sl The chirality group $X({\cal H})$ of each orientably regular
hypermap $\cal H$ is isomorphic to a normal
subgroup of the monodromy group $\Mon({\cal H})$.}
\vskip5pt
\newcounter{cgsg}
\setcounter{cgsg}{\theteocounter}
\addtocounter{teocounter}{1}
\noindent {\sl Proof.}
$X({\cal H})\cong H^{\Delta}/H\leq\Delta^+/H\cong\Mon({\cal H})$.
By Theorem~\theskfc,\ $X$ is normal in $\Mon(\cal H)$.
\hspace{\fill} $\Box$
\vskip5pt

Since the number of darts in an orientably regular hypermap
coincides with the order of the monodromy group, Proposition \thecgsg\ and
Lagrange's Theorem
now imply:

\vskip7pt
\noindent {\bf Corollary \theteocounter\,\,}
{\sl The chirality index of any orientably regular hypermap divides
the number of darts.}
\newcounter{IO}
\setcounter{IO}{\theteocounter}
\addtocounter{teocounter}{1}

\vskip7pt
The following example shows that the chirality
index of a hypermap can be arbitrarily large.
\vskip7pt

\noindent {\bf Example 1} (see \cite{CHIF}) Consider the
metacyclic group
$$G=\langle\, a,b \mid a^n=1,\, b^m=a^s,\,
bab^{-1}=a^r\,\rangle$$
of order $mn$, where $rs\equiv s$ (mod~$n$) and $r^m\equiv
1$ (mod~$n$) for some $m>2$. If we regard $a$ and $b$
as the canonical generators
for an orientably regular hypermap ${\cal H}=(G,a,b)$ with
monodromy group $G$, then the chirality group $X({\cal H})$
can be interpreted as the smallest normal subgroup $N$ of $G$
such that the assignment $a\mapsto a^{-1}$ and $b\mapsto b^{-1}$
induces an automorphism of $G/N$. We
obtain this quotient from $G$ by adding the extra relations
formed from those of $G$ by replacing $a$ and $b$ with their
inverses. In this case, it is sufficient to add
$b^{-1}a^{-1}b=a^{-r}$, or equivalently, $b^{-1}ab=a^r$, so that
$a=b(b^{-1}ab)b^{-1}=(a^r)^r=a^{r^2}$ in $G/N$. Thus
$a^{r^2-1}=1$ in $G/N$, so in $G$ it follows that $K=\langle
a^{r^2-1}\rangle$ is a
subgroup of $N$. On the other hand, it is easy to see that $K$ is a normal
subgroup of $G$ such that $G/K$ is invariant under replacement
of the generators with their inverses. By minimality  $N=K$, so $X({\cal
H})\cong \langle
a^{r^2-1}\rangle$ and $\kappa({\cal H})=n/\gcd(n,r^2-1)$; since $m>2$, this
can be arbitrarily
large.

\section{Normal Subgroups of $\Delta$ and $\Delta^+$}

In order to study chirality groups further we need some
elementary results about which normal subgroups of $\Delta^+$ are
also normal in $\Delta$.

The following result is obvious.

\vskip8pt
\noindent {\bf Lemma \theteocounter\,\,} {\em Let $N$ be a normal subgroup
of $\Delta^+$,
and let $G=\Delta^+/N$. Then the following are equivalent:
\begin{description}
\item[(i)] $N$ is normal in $\Delta$;
\item[(ii)] $N^r=N$;
\item[(iii)] $G$ has an automorphism inverting both of its canonical
generators;
\item[(iv)] the hypermap corresponding to $N$ is reflexible.
\end{description}
}
\newcounter{AA}
\setcounter{AA}{\theteocounter}
\addtocounter{teocounter}{1}

Let us call a generating pair $x,y$ for a group $G$ {\em
symmetric\/} if there is an automorphism of $G$ inverting both
$x$ and $y$, and {\em asymmetric\/} otherwise. Let us call a
$2$-generator group $G$ {\em strongly symmetric\/} if all its
generating pairs are symmetric. By Lemma \theAA, this condition means
that every normal subgroup of $\Delta^+$ with quotient group $G$
is normal in $\Delta$, or equivalently, every orientably regular
hypermap with monodromy group $G$ is regular.

\vskip8pt
\noindent {\bf Proposition \theteocounter\,\,} {\sl The following groups
are all strongly
symmetric: $2$-generator abelian groups, dihedral groups,
$PSL_2(q)$ for any prime power $q$, and the symmetric group $S_n$
for $n\leq 5$.}
\vskip5pt
\newcounter{BB}
\setcounter{BB}{\theteocounter}
\addtocounter{teocounter}{1}
\noindent{\sl Proof.} Every abelian group has an automorphism
inverting all its elements. If two elements generate a dihedral
group, they are either two reflections or a reflection and a
rotation; in the first case they are both inverted by the identity
automorphism, and in the second case they are inverted by
conjugation by the reflection. Macbeath \cite{Mac} showed that if $x,
y$ and $x', y'$ are two generating pairs for $PSL_2(q)$ with
${\rm tr}\,x={\rm tr}\,x'$, ${\rm tr}\,y={\rm tr}\,y'$ and ${\rm
tr}\,xy={\rm tr}\,x'y'$, then there is an automorphism sending
$x$ to $x'$ and $y$ to $y'$; now ${\rm tr}\,x={\rm tr}\,x^{-1}$,
${\rm tr}\,y={\rm tr}\,y^{-1}$ and ${\rm tr}\,xy={\rm
tr}\,yx={\rm tr}\,(yx)^{-1}={\rm tr}\,x^{-1}y^{-1}$, so $x, y$
form a symmetric pair (see \cite{Sin} for this and related arguments).
The small symmetric groups are dealt with by a routine inspection
of their generating pairs.
\hspace{\fill} $\Box$
\vskip5pt

Note that this shows that the alternating groups $A_n$ are
strongly symmetric for $n\leq 6$ since $A_1$ and $A_2$ are
trivial, $A_3\cong C_3$, $A_4\cong PSL_2(3)$, $A_5\cong PSL_2(4)$
and $A_6\cong PSL_2(9)$.

\vskip8pt

An immediate consequence of Proposition \theBB\, is the following:

\vskip8pt
\noindent {\bf Corollary \theteocounter\,\,} {\sl If $N$ is a normal
subgroup of $\Delta^+$
such that $\Delta^+/N$ is abelian or dihedral or isomorphic to
$PSL_2(q)$ or $S_n$ for some $n\leq 5$, then $N$ is normal in
$\Delta$.}
\newcounter{CC}
\setcounter{CC}{\theteocounter}
\addtocounter{teocounter}{1}
\vskip8pt

\noindent{\bf Example 2.} As an application of Corollary \theCC, the
following simple example shows that every elementary abelian
group of prime power order $q=p^e>4$ can arise as a chirality
group. (We will extend this result later to all finite abelian
groups, but for this we need a more complicated construction.)
Let $G$ be the 1-dimensional affine group $AGL_1(q)$ over the
field $F_q$ of order $q$, consisting of the transformations
$t\mapsto at+b$ where $a,b\in F_q$ and $a\neq 0$. This is a split
extension of the translation group $B=\{t\mapsto t+b\mid b\in
F_q\}$, an elementary abelian group of order $q$ isomorphic to
the additive group of $F_q$, by the group $A=\{t\mapsto at\mid
a\in F_q^*\}$, isomorphic to the multiplicative group
$F_q^*=F_q\setminus\{0\}$ of $F_q$. Let $x$ and $y$ be the
elements $t\mapsto ct$ and $t\mapsto t+1$ of $G$, where $c$
generates the (cyclic) group $F_q^*$. Then $x$ and $y$ generate
$G$, so the epimorphism $\theta:\Delta^+\to G,\; R\mapsto x,
L\mapsto y$ realises $G$ as the monodromy group of an orientably
regular hypermap $\cal H$, with hypermap subgroup $H={\rm
ker}\,\theta\leq\Delta^+$.

If $\Gamma$ denotes the Galois group of $F_q$ (over its prime
field $F_p$), then ${\rm Aut}\,G$ can be identified with $A\Gamma L_1(q)$,
the group of transformations $t\mapsto at^{\gamma}+b$
where $a\neq 0$ and $\gamma\in\Gamma$; this contains $G$ as a
normal subgroup, and induces automorphisms of $G$ by acting by
conjugation. It follows that $x$ is inverted by an automorphism
of $G$ if and only if $c$ is equivalent to $c^{-1}$ under
$\Gamma$. Now $\Gamma$ is a cyclic group of order $e$, generated
by the Frobenius automorphism $t\mapsto t^p$ of $F_q$, so the
images of $c$ under $\Gamma$ are the powers $c^{p^i}$ for
$i=1,\ldots, e$. Thus $c$ is equivalent to $c^{-1}$ if and only
if $c^{p^i+1}=1$ for some $i$, or equivalently (since $c$ has
order $p^e-1$) if $p^e-1$ divides $p^i+1$, which is impossible if
$p^e>4$. Thus $\cal H$ is chiral for all $q>4$, so $H^{\Delta}>H$.

If $N$ denotes $\theta^{-1}(B)$ then $\Delta^+/N\cong G/B\cong
A$, which is abelian, so $N$ is normal in $\Delta$ by
Corollary \theCC, and hence $H^{\Delta}\leq N$. Now $B$ is a minimal
normal subgroup of $G$, so $N/H$ is a minimal normal subgroup of
$\Delta^+/H$. Since $N\geq H^{\Delta}>H$ it follows that
$H^{\Delta}=N$, so $\cal H$ has chirality group $X({\cal H})\cong
N/H\cong B$.
\vskip5pt

\section{Totally Chiral Hypermaps}

The most extreme type of chirality, and often the easiest to
study, occurs when an orientably regular hypermap $\cal H$ is {\em totally
chiral\/}, meaning that the chirality group of $\cal H$ coincides with its
monodromy
group. If $H$ is the map subgroup of $\Delta^+$ corresponding to $\cal H$,
then $\cal H$ is totally
chiral if and only if $H^{\Delta}={\Delta}^+$, so a totally chiral hypermap
is one which covers no
nontrivial reflexible hypermaps.

In this section we will study some examples of this phenomenon,
before considering more general forms of chirality in the next
section.

\vskip8pt

A {\em perfect\/} group is one with no nontrivial abelian
epimorphic images.

\vskip8pt
\noindent {\bf Corollary \theteocounter\,\,} {\sl The monodromy group of a
totally chiral
hypermap cannot have any nontrivial strongly symmetric group as
an epimorphic image; in particular, it must be perfect.}
\vskip5pt
\newcounter{DD}
\setcounter{DD}{\theteocounter}
\addtocounter{teocounter}{1}
\noindent{\sl Proof.}
Let $\cal H$ be a totally chiral hypermap with monodromy group
$\Mon({\cal H})={\Delta}^+/H$ for some $H\unlhd{\Delta}^+$.
Assume, for a contradiction, that $\Mon({\cal H})$ has a nontrivial
strongly symmetric epimorphic image $G$. Then there exists a
subgroup $N$ of ${\Delta}^+$ such that $H\leq N < {\Delta}$ and
$G={\Delta}^+/N$. Since ${\Delta}^+/N$ is strongly symmetric,
we have $N\unlhd\Delta$ and hence $NN^r=N$. Therefore
$${\Delta}^+=HH^r\leq NN^r=N<{\Delta}^+$$
which is a contradiction. In particular, since abelian groups are strongly
symmetric, they cannot occur as nontrivial epimorphic images.
\hspace{\fill} $\Box$
\vskip5pt

The most obvious examples of perfect groups are the nonabelian
simple groups, so it makes sense to inspect these for examples of
totally chiral hypermaps.

\vskip8pt
\noindent {\bf Lemma \theteocounter\,\,} {\em If $X$ is a simple group with
an asymmetric
generating pair, then $X$ is the monodromy group of a totally
chiral hypermap.}
\vskip5pt
\newcounter{EE}
\setcounter{EE}{\theteocounter}
\addtocounter{teocounter}{1}
\noindent{\sl Proof.} Let $x$ and $y$ be an asymmetric generating
pair for $X$, and let $H$ be the kernel of the epimorphism
$\Delta^+\to X$ given by $R\mapsto x, L\mapsto y$, so
$X\cong\Delta^+/H$ is the monodromy group of the
orientably regular hypermap $\cal H$ corresponding to $H$. By Lemma \theAA\, we
have $H^r\neq H$, so $H^{\Delta}$ is a normal subgroup of
$\Delta^+$ properly containing $H$. The simplicity of $X$ gives
$H^{\Delta}=\Delta^+$, so $\cal H$ is totally chiral.
\hspace{\fill} $\Box$
\vskip5pt

Some simple groups $X$ have an element $x$ which is not inverted
by any automorphism; if there is a second element $y$ such that
$x$ and $y$ generate $X$, then Lemma \theEE\, implies that $X$ is the
monodromy group of a totally chiral hypermap. Here we give two
families of groups which illustrate this principle.

\vskip8pt

The Ree groups $Re(3^f)=\/^2\negthinspace G_2(3^f)$ and the
Suzuki groups $Sz(2^f)=\/^2\negthinspace B_2(2^f)$ are defined
for all odd $f\geq 1$, and are simple for $f>1$. They were first
described in \cite{Re1, Re2} and \cite{Suz} respectively, and a good
account of their properties can also be found in Chapter XI of
\cite{HB}.

\vskip8pt \noindent {\bf Theorem \theteocounter\,\,} {\em For
each odd $f>1$ there are totally chiral hypermaps with monodromy
groups $Re(3^f)$ and $Sz(2^f)$.} \vskip5pt
\newcounter{FF}
\setcounter{FF}{\theteocounter}
\addtocounter{teocounter}{1}
\noindent{\sl Proof.} The Suzuki group $X=Re(3^f)$, for odd
$f>1$, is a simple epimorphic image of the triangle group
$\Delta(2,3,7)$, and hence of $\Delta^+$, with generators $x$ and
$y$ of orders $2$ and $3$. As shown in \cite{Re2}, $y$ is not inverted
by any automorphism of $X$, so Lemma \theEE\, applies. Similarly
$Sz(2^f)$ is a simple epimorphic image of $\Delta(2,4,5)$ for odd
$f>1$, and the generator of order $4$ is not inverted by any
automorphism \cite{Suz}. (See \cite[Section 6]{Jon}  and \cite[Section 6]{JS}
for further details and for descriptions of the chiral hypermaps
associated with these groups.)
\hspace{\fill} $\Box$
\vskip5pt
Other simple examples of this phenomenon include the Mathieu
groups $M_p$ for $p=11$ and $23$: in each case an element of
order $p$ is a member of a generating pair, but is not inverted
by any automorphism (see \cite{CCNPW} for details).

\vskip8pt

A more common phenomenon among simple groups is that both members
$x$ and $y$ of a generating pair may be inverted by
automorphisms, but no single automorphism inverts them both. To
illustrate this we consider the alternating groups $A_n$, which
are simple for all $n\geq 5$; here every element, having the same
cycle structure as its inverse, is inverted by conjugation in
$S_n$.

\vskip8pt

\vskip8pt
\noindent {\bf Theorem \theteocounter\,\,} {\em If $n\geq 7$ there is a
totally chiral
hypermap with monodromy group $A_n$.}
\vskip5pt
\newcounter{GG}
\setcounter{GG}{\theteocounter}
\addtocounter{teocounter}{1}
\noindent{\sl Proof.} If $n$ is odd, let $X$ be the subgroup of
the symmetric group $S_n$ generated by the permutations $x=(1\,
2\,\ldots\, n)$ and $y=(1\, 2\, 4)$. These are both even, so
$X\leq A_n$. Any $x$-invariant relation must be congruence
mod~$(m)$ for some $m$ dividing $n$, and $y$ preserves this only
for $m=1$ and $n$, so $X$ is primitive. A primitive group
containing a 3-cycle must contain $A_n$ \cite[Theorem 13.3]{Wie}, so
$X=A_n$. For $n\geq 7$ one can identify ${\rm Aut\,}A_n$ with
$S_n$, acting by conjugation on its normal subgroup $A_n$; the
only permutations inverting $x$ are the reflections in the
obvious dihedral group $D_n$ containing $x$, and these do not
invert $y$, so Lemma \theEE\, gives the required result.

If $n$ is even, one can take $x=(1\, 2\,\ldots\, n-1)$ and
$y=(1\, n)(2\, 3)$ in $A_n$. The group $X=\langle x,y\rangle$ is
clearly transitive; since $x$ fixes $n$ and has a single cycle on
the remaining points, $X$ is doubly transitive and hence
primitive. Now $X$ contains a $5$-cycle
$[x,y]=x^{-1}y^{-1}xy=(1\, n\, 3\, 4\, 2)$, and a theorem of
Jordan \cite[Theorem 13.9]{Wie} states that a primitive group of
degree $n$, containing a $p$-cycle for some prime $p\leq n-3$,
must contain $A_n$. It follows that $X=A_n$ provided $n\geq 8$.
Again there is no permutation in $S_n$ which inverts $x$ and $y$,
so Lemma \theEE\, completes the proof.
\hspace{\fill} $\Box$
\vskip5pt

By the comment after Proposition \theBB, this result does not extend
to $A_n$ for $n\leq 6$.

\vskip8pt

For a similar class of examples, we use the simple groups
$PSL_d(q)$ and their covering groups $SL_d(q)$. Note that the
case $d=2$ is excluded by Proposition \theBB\, and Corollary \theDD, and
similarly $PSL_3(2)=SL_3(2)$ is excluded since it is isomorphic
to $PSL_2(7)$.

\vskip8pt

Each $a\in GL_d(q)$ induces an inner automorphism
$\iota_a:g\mapsto g^a=a^{-1}ga$ of $GL_d(q)$, each
$\gamma\in\Gamma={\rm Gal}\,F_q$ induces an automorphism by
acting on matrix entries, and the adjoint mapping $g\mapsto
g^*=(g^T)^{-1}$ is also an automorphism of $GL_d(q)$ (the {\em
graph automorphism\/} in Lie algebra terminology). Every
automorphism of $GL_d(q)$ is a composition of these, having the
form $\alpha_{a,\gamma}:g\mapsto (g^a)^{\gamma}$ or
$\beta_{a,\gamma}:g\mapsto ((g^*)^a)^{\gamma}$. These
automorphisms all induce automorphisms of $PGL_d(q),\, SL_d(q)$ and
$PSL_d(q)$, and conversely all automorphisms of these groups
arise in this way. Since every matrix is similar (that is,
conjugate) to its transpose, every element of these groups is
inverted by some automorphism. Nevertheless, it is usually
possible to find an asymmetric generating pair in these groups.

\vskip8pt
\noindent {\bf Theorem \theteocounter\,\,} {\em If $d\geq 3$ and if the
prime power $q$ is
sufficiently large (in terms of $d$) then there are totally
chiral hypermaps with monodromy groups $SL_d(q)$ and
$PSL_d(q)$.}
\vskip5pt
\newcounter{HH}
\setcounter{HH}{\theteocounter}
\addtocounter{teocounter}{1}
\noindent{\sl Proof.} Let $x$ be a diagonal matrix in
$X=SL_d(q)$, with spectrum $\Lambda\subset F_q$ consisting of $d$
distinct eigenvalues, so the centraliser of $x$ in $GL_d(q)$
consists of the diagonal matrices. We will show that if $q$ is
sufficiently large then one can choose $\Lambda$ so that $x$ is a
member of an asymmetric generating pair $x, y$ for $X$, with a
similar result for their images $\overline x,\overline y$ in
$\overline X=PSL_d(q)$.

The automorphisms $\iota_a$ preserve eigenvalues, the
field-automorphisms act naturally on them, and the adjoint
automorphism inverts them. The image $\overline x$ of $x$ in
$\overline X$ is represented by the matrices $\omega x$ where
$\omega^d=1$, and these have spectrum $\omega\Lambda$. For any
fixed $d\geq 3$, if $q$ is sufficiently large one can choose
$\Lambda$ to be distinct from the sets
$\omega(\Lambda^{-1})^{\gamma}$ where $\omega^d=1$ and
$\gamma\in\Gamma$, and also from the sets
$\omega\Lambda^{\gamma}$ except where $\omega=1$ and $\gamma=1$.
(For example, one can take
$\Lambda=\{\zeta,\zeta^2,\ldots,\zeta^{d-1},\zeta^{-d(d-1)/2}\}$
for sufficiently large $q$, where $\zeta$ generates $F_q^*$,
noting that inversion and multiplication by $\omega$ preserve the
multiset of ratios $\lambda/\mu$ for $\lambda,\mu\in\Lambda$; no
such choice is possible for $d=2$ since the condition $\det x=1$
forces $\Lambda=\Lambda^{-1}$.) Our conditions on $\Lambda$ then
ensure that the only automorphisms inverting $\overline x$ (and
hence the only automorphisms inverting $x$) are the automorphisms
$\beta_{a,1}: g\mapsto (g^*)^a$ where $a$ commutes with $x$ and
is therefore a diagonal matrix. A simple calculation shows that
the proportion of elements $y\in X$ (or $\overline y\in\overline
X$) inverted by any such $\beta_{a, 1}$ approaches $0$ as
$q\to\infty$, whereas recent results on random generation of
simple groups \cite{GKS} imply that the proportion of $\overline y$
such that $\langle \overline x,\overline y\rangle=\overline X$
approaches $1$. This means that for all sufficiently large $q$,
$\overline X$ has an asymmetric generating pair $\overline
x,\overline y$, so Lemma \theEE\, shows that $\overline X$ is the
monodromy group of a totally chiral hypermap, corresponding to
the kernel $H$ of the obvious epimorphism $\Delta^+\to\overline
X$.

Clearly $x,y$ form an asymmetric pair, since any automorphism of
$X$ inverting them would induce an automorphism of $\overline X$
inverting $\overline x$ and $\overline y$; we now show that they
generate $X$. The subgroup $S=\langle x,y\rangle$ maps onto
$\overline X$, so $SZ=X$ where $Z$ is the centre of $X$, and
hence $S$ is normal in $X$ (since it commutes with $Z$); now
$X/S=SZ/S\cong Z/(S\cap Z)$, which is abelian, and $X$ is
perfect, so $S=X$ as required. This gives an epimorphism
$\Delta^+\to X$ with kernel $K$ contained in $H$. Now
$K^{\Delta}$ is a normal subgroup of $\Delta^+$, properly
containing $K$, so by the normal structure of $X=SL_d(q)$ we have
$K^{\Delta}\leq H$ or $K^{\Delta}=\Delta^+$. The first case
implies that $H/K^{\Delta}$ is the centre of
$\Delta^+/K^{\Delta}$, so $H$ is normal in $\Delta$, which is
false. Hence $K^{\Delta}=\Delta^+$, so the hypermap corresponding
to $K$ is totally chiral, with monodromy group $X$.
\hspace{\fill} $\Box$
\vskip5pt

No doubt, a more careful examination of the families of finite
simple groups would reveal further examples of totally chiral
hypermaps.

\section{Chirality Groups and Direct Products}

In this section we define a form of direct product of hypermaps, and we
investigate how the chirality group of the product can be
expressed in terms of the chirality groups of the factors. As an
application of our general results we show that every finite
abelian group is the chirality group of some hypermap.

Let $\cal H$ and $\cal K$ be orientably regular hypermaps with hypermap
subgroups $H,K\leq\Delta^+$. We define the {\em least common cover\/} ${\cal
H}\vee{\cal K}$ and the {\em greatest common quotient\/} ${\cal
H}\wedge{\cal K}$ of $\cal H$ and
$\cal K$ to be the orientably regular hypermaps with hypermap subgroups
$H\cap K$ and $\langle
H,K\rangle=HK$ respectively. For example, we have ${\cal H}_{\Delta}={\cal
H}\vee{\cal H}^r$ and
${\cal H}^{\Delta}={\cal H}\wedge{\cal H}^r$. Clearly, every common cover
of $\cal H$ and $\cal
K$ also covers ${\cal H}\vee{\cal K}$, and any
common quotient of them is also a quotient of ${\cal H}\wedge{\cal K}$. The
coverings ${\cal H}\vee{\cal K}\to{\cal H}$ and
${\cal K}\to{\cal H}\wedge{\cal K}$ are
regular with covering transformation group
$H/(H\cap K)\cong HK/K$, and the coverings ${\cal H}\vee{\cal K}\to{\cal
K}$ and ${\cal H}\to{\cal
H}\wedge{\cal K}$ are regular with covering transformation
group $K/(H\cap K)\cong HK/H$.

We can write ${\cal H}=(D_1,R_1,L_1)$ and ${\cal K}=(D_2,R_2,L_2)$ for sets
$D_i$ of darts, so
that the actions $\rho\mapsto R_i,\,\lambda\mapsto L_i$ of $\Delta^+$ on
$D_1$ and $D_2$ induce a
product action
$\rho\mapsto R,\,\lambda\mapsto L$ of
$\Delta^+$ on the set $D=D_1\times D_2$, given by $R(x,y)=(R_1x,R_2y)$ and
$L(x,y)=(L_1x,L_2y)$.
If this action is transitive then it determines an orientably regular hypermap
${\cal H}\times{\cal K}=(D,R,L)$ called the {\em oriented direct product\/}
of ${\cal H}$ and
${\cal K}$, with hypermap subgroup $H\cap K$. The following straightforward
lemma tells us when
this construction is possible (see \cite{Bre}):

\medskip
\vskip8pt
\noindent {\bf Lemma \theteocounter\,\,} {\em If $\cal H$ and $\cal K$ are
orientably regular
hypermaps, then the following three conditions are equivalent:}
\begin{description}
\item[(a)] {\em $\Delta^+$ acts transitively on $D$;}

\item[(b)] {\em ${\cal H}\wedge{\cal K}$ is the trivial orientable hypermap,
with one dart;}

\item[(c)] $HK=\Delta^+$.
\end{description}
\vskip5pt
\newcounter{product}
\setcounter{product}{\theteocounter}
\addtocounter{teocounter}{1}

When these conditions are satisfied, we say that $\cal H$ and
$\cal K$ are {\em orientably orthogonal\/}, and write ${\cal H}\perp{\cal
K}$; then ${\cal
H}\times{\cal K}$ exists, and is isomorphic to ${\cal H}\vee{\cal K}$, with
monodromy group
$\Mon({\cal H}\times{\cal K})=\Mon({\cal H})\times\Mon({\cal K})$. For example,
$\cal H$ is totally chiral if and only if ${\cal H}\perp{\cal H}^r$, in
which case ${\cal
H}_{\Delta}={\cal H}\times{\cal H}^r$ with $\Mon({\cal
H}_{\Delta})=\Mon({\cal H})\times\Mon({\cal
H}^r)\cong (\Mon({\cal H}))^2$.

\medskip

\vskip8pt
\noindent {\bf Proposition \theteocounter\,\,} {\em Let $\cal H$ and $\cal
K$ be orientably
regular hypermaps, with hypermap subgroups $H$ and $K$, such that $\cal K$
is totally chiral and
covers $\cal H$. Then the product ${\cal L}={\cal K}\times{\cal H}^r$ is an
orientably regular
hypermap with chirality group $X({\cal L})\cong H/K$}
\vskip5pt
\newcounter{factor}
\setcounter{factor}{\theteocounter}
\addtocounter{teocounter}{1}
{\sl Proof.} Since $\cal K$ is totally chiral and $K\le H$, we have
$KH^r\ge KK^r=\Delta^+$, so
$KH^r=\Delta^+$ and hence ${\cal K}\perp {\cal H}^r$. Consequently, the product
${\cal L}={\cal K}\times{\cal H}^r$ is a well-defined orientably regular
hypermap with
hypermap subgroup $L=K\cap H^r$. Using $\core{L}{\Delta}=K\cap L^r$, the
third isomorphism theorem,
and $KL^r=H$, we have
$X({\cal L})\cong L^r/L_{\Delta}=L^r/(K\cap L^r)\cong KL^r/K=H/K$.
\hspace{\fill} $\Box$
\vskip5pt

We define two groups $A$ and $B$ to be {\em mutually orthogonal\/},
and write $A\perp B$, if they have no non-trivial common epimorphic images.

\vskip8pt
\noindent {\bf Lemma \theteocounter\,\,}
{\em Let $\cal H$ and $\cal K$ be orientably regular hypermaps.}
\begin{description}
\item[(a)] {\em If} $\Mon({\cal H})\perp\Mon({\cal K})$, {\em then ${\cal
H}\perp{\cal K}$.}
\item[\rm(b)] {\em If, in addition, $\cal H$ and $\cal K$ are totally chiral,
so is ${\cal H}\vee{\cal K}$.}
\end{description}
\newcounter{comim}
\setcounter{comim}{\theteocounter}
\vskip5pt
\addtocounter{teocounter}{1}
\noindent{\sl Proof.} (a) Let $H$ and $K$ be the hypermap subgroups of
$\Delta^+$
corresponding to $\cal H$ and $\cal K$. If $HK\neq\Delta^+$, then
$\Delta^+/HK$ is a
non-trivial epimorphic image of $\Delta^+/H\cong\Mon({\cal H})$ and of
$\Delta^+/K\cong\Mon({\cal
K})$, against our assumption.

\medskip

\noindent(b) Since ${\cal H}\perp{\cal H}^r$, we have ${\cal
H}_{\Delta}={\cal H}\times{\cal H}^r$,
and similarly for ${\cal K}_{\Delta}$. The groups $\Mon({\cal
H}_{\Delta})=\Mon({\cal H})
\times\Mon ({\cal H}^r)\cong(\Mon({\cal H}))^2$ and $\Mon({\cal
K}_{\Delta})=\Mon({\cal K})
\times\Mon({\cal K}^r)\cong(\Mon({\cal K}))^2$ are mutually orthogonal
(since $\Mon({\cal H})$ and $\Mon({\cal K})$ are), so ${\cal
H}_{\Delta}\perp{\cal K}_{\Delta}$
by (a). Now $(H\cap K)_{\Delta} = H_{\Delta}\cap K_{\Delta}$, so
$$\Mon({\cal H}\vee{\cal K})_{\Delta}=\Mon({\cal
H}_{\Delta})\times\Mon({\cal K}_{\Delta})$$
$$=(\Mon({\cal H})\times\Mon({\cal H}^r))\times(\Mon({\cal
K})\times\Mon({\cal K}^r))$$
$$=(\Mon({\cal H})\times\Mon({\cal K}))\times(\Mon({\cal
H}^r)\times\Mon({\cal K}^r))$$
$$=\Mon({\cal H}\vee{\cal K})\times\Mon(({\cal H}\vee{\cal K})^r)$$
as required.
\hspace{\fill} $\Box$
\vskip5pt
The converse of (a) is false: for instance, for each prime $p$ there are
$p+1$ mutually orthogonal orientably regular hypermaps $\cal H$ with
$\Mon({\cal H})\cong C_p$.
\vskip8pt
\noindent {\bf Corollary \theteocounter\,\,}
{\em Let ${\cal H}_1,\ldots,{\cal H}_k$ be orientably regular hypermaps.}
\begin{description}
\item[(a)] {\em If the groups $\Mon({\cal H}_i)$ are pairwise orthogonal, then
\newline ${\cal H}_1\vee\cdots\vee{\cal H}_k={\cal
H}_1\times\cdots\times{\cal H}_k$.}
\item[(b)] {\em If, in addition, each hypermap ${\cal H}_i$ is totally
chiral, then so is ${\cal H}_1\vee\cdots\vee{\cal H}_k$.}
\end{description}
\newcounter{ort}
\setcounter{ort}{\theteocounter}
\vskip5pt
\addtocounter{teocounter}{1}

\noindent{\sl Proof.} This follows by induction on $k$: since
$\Mon({\cal H}_i)\perp \Mon({\cal H}_k)$ for all $i<k$ we have
$\Mon({\cal H}_1)\times\cdots\times\Mon({\cal H}_{k-1})\perp
\Mon({\cal H}_k)$, so the previous Lemma applies.
\hspace{\fill} $\Box$
\vskip5pt

We are now ready to prove that every finite abelian group appears as the
chirality group
of some hypermap. First we deal with cyclic groups:
\vskip8pt
\noindent {\bf Proposition \theteocounter\,\,} {\em For each integer $n\geq
1$ there is an
orientably regular hypermap with chirality group $C_n$.}
\vskip5pt
\newcounter{II}
\setcounter{II}{\theteocounter}
\addtocounter{teocounter}{1}
\noindent{\sl Proof.} We have seen in the proof of Theorem \theHH\, that if
$d\geq 3$ and $q$ is sufficiently large, then there are
totally chiral hypermaps ${\cal K}$ and ${\cal H}$ with
normal hypermap subgroups $K\leq H$ in $\Delta^+$ satisfying $\Delta^+/K\cong
X=SL_d(q)$, $\Delta^+/H\cong\overline X=PSL_d(q)$, and $H/K\cong
Z=Z(X)$. By Proposition~\thefactor, $Z\cong X({\cal L})$ where
${\cal L}={\cal K}\times {\cal H}^r$.
Now $Z$ is a cyclic group of order $\gcd(q-1,d)$, and we can choose
$d\;(\geq 3)$ and $q$ (sufficiently large) so that this takes any
prescribed value $n\geq 1$: if $n\geq 3$ we can take $d=n$ and
$q\equiv 1$ (mod~$n$) (recall that Euler's Theorem gives
$p^{\phi(n)}\equiv 1$ (mod~$n$) for each prime $p$ not dividing
$n$), and if $n=1$ or $2$ we can take $d=4$ and $q\equiv 0$ or
$3$ (mod~$4$), respectively.

\hspace{\fill} $\Box$

\vskip8pt

We now extend this method of proof to all abelian groups:

\vskip8pt
\noindent {\bf Theorem \theteocounter\,\,} {\em Every finite abelian group
is the chirality
group of some orientably regular hypermap.}
\vskip5pt
\newcounter{JJ}
\setcounter{JJ}{\theteocounter}
\addtocounter{teocounter}{1}
\noindent{\sl Proof.} Every finite abelian group can be written in the form
$A=C_{n_1}\times\cdots\times C_{n_k}$ for some positive integers
$n_1,\ldots, n_k$. By the proof of Theorem~\theHH, we can choose normal
subgroups
$K_i\leq H_i$ of $\Delta^+$ for $i=1,\ldots, k$ such that
$\Delta^+/K_i\cong SL_{d_i}(q_i),\,\Delta^+/H_i\cong PSL_{d_i}(q_i)$ and
$H_i/K_i=Z(\Delta^+/K_i)\cong C_{n_i}$ with $d_i\geq 3$, and
the corresponding orientably regular hypermaps ${\cal K}_i$ and ${\cal
H}_i$ are totally
chiral; moreover, if $n_i=n_j$ for some $i\neq j$ then we can take $q_i\neq
q_j$. Now every proper
normal subgroup of $SL_{d_i}(q_i)$ is central, so every non-trivial
epimorphic image of
$SL_{d_i}(q_i)$ maps onto $PSL_{d_i}(q_i)$; for distinct pairs $d_i\;(\geq
3)$ and $q_i$ the simple
groups $PSL_{d_i}(q_i)$ are non-isomorphic, so they are mutually
orthogonal, and hence so are the
groups $SL_{d_i}(q_i)$.

By Corollary \theort, the hypermaps ${\cal H}={\cal H}_1\vee\cdots\vee{\cal
H}_k={\cal H}_1\times\cdots\times{\cal H}_k$ and ${\cal K}={\cal
K}_1\vee\cdots\vee{\cal K}_k={\cal K}_1\times\cdots\times{\cal K}_k$ are
both totally chiral,
corresponding to normal subgroups $H=H_1\cap\cdots\cap H_k$ and
$K=K_1\cap\cdots\cap K_k$ of
$\Delta^+$ with $H/K\cong(H_1/K_1)\times\cdots\times(H_k/K_k)\cong A$. By
Proposition
\thefactor\, the hypermap ${\cal L}={\cal K}\times{\cal H}^r$ corresponding
to $L=K\cap H^r$ has
$X({\cal L})\cong H/K\cong A$.
\hspace{\fill} $\Box$
\vskip5pt
\section{Nonabelian chirality groups}
\vskip5pt
The preceding results make it tempting to conjecture that every
finite group is the chirality group of some orientably regular hypermap.
The next result shows that
this is false. For any group $S$, let $Z, I$ and $A$ denote the
centre $Z(S)$, the inner automorphism group ${\rm Inn}\,S\cong
S/Z$, and the automorphism group ${\rm Aut}\,S$. The action of
$S$ by conjugation on itself induces a homomorphism $S\to I\leq
A$ with kernel $Z$; for each subgroup $T\leq S$ let $\overline T$
denote the image $TZ/Z$ of $T$ in $I$. A subgroup $T\leq S$ is
{\it characteristic\/} in $S$ if it is invariant under $A$; this
implies that $\overline T$ is a normal subgroup of $A$.

\vskip8pt
\noindent {\bf Theorem \theteocounter\,\,} {\sl If a group $S$ has a
characteristic subgroup
$T$ such that $\overline T<\overline S$ and $A/\overline T$ is
abelian, then $S$ cannot be a chirality group.}
\vskip5pt
\newcounter{KK}
\setcounter{KK}{\theteocounter}
\addtocounter{teocounter}{1}

Before proving this result, let us examine some of its
consequences. Recall that a group $S$ is {\it complete\/} if $S$
has trivial centre and all its automorphisms are inner, so that
$S\cong {\rm Aut}\,S$.

\vskip8pt

\vskip8pt
\noindent {\bf Corollary \theteocounter\,\,} {\sl If $S$ is a complete
group and $S$ is not
perfect, then $S$ cannot be a chirality group.}
\vskip5pt
\newcounter{LL}
\setcounter{LL}{\theteocounter}
\addtocounter{teocounter}{1}
\noindent{\sl Proof.} Let $T$ be the commutator subgroup $S'$, so
$T<S$ since $S$ is not perfect. Since $Z=1$ we have $\overline
S=S$ and $\overline T=T$, so $\overline T<\overline S$. Since
$A=S$ we have $A/\overline T=S/S'$, which is abelian, so
Theorem \theKK\, gives the result.
\hspace{\fill} $\Box$
\vskip5pt

\vskip8pt
\noindent {\bf Corollary \theteocounter\,\,} {\sl The symmetric group $S_n$
is a chirality
group if and only if $n\leq 2$.}
\vskip5pt
\newcounter{MM}
\setcounter{MM}{\theteocounter}
\addtocounter{teocounter}{1}
\noindent{\sl Proof.} Let $S=S_n$ for any $n>2$, and take $T$ to
be the alternating group $A_n$, a proper subgroup of $S$ which is
characteristic since $T=S'$. We have $Z=1$, and if $n\neq 6$ then
$A=S$, so Corollary \theLL\, implies that $S$ cannot be a chirality
group. If $n=6$ then $|A:S|=2$, so $A/\overline T=A/T$ is abelian
(of order $4$) and Theorem \theKK\, gives the required result. However,
$S_n\cong C_n$ for $n=1$ and $2$, and these are chirality groups
by Theorem \theII.
\hspace{\fill} $\Box$
\vskip5pt

More generally, the automorphism group $S$ of a nonabelian finite
simple group $F$ is always complete. In many cases $F$ has outer
automorphisms, so $S>F$; by the Schreier Conjecture $S/F$ is
solvable, so $S$ is not perfect and hence cannot be a chirality
group.

\vskip8pt
\noindent {\bf Corollary \theteocounter\,\,} {\sl The dihedral group $D_n$
is a chirality
group if and only if $n\leq 2$.}
\vskip5pt
\newcounter{NN}
\setcounter{NN}{\theteocounter}
\addtocounter{teocounter}{1}
\noindent{\sl Proof.} Let $S$ be the dihedral group $D_n=\langle
a,b\mid a^n=b^2=(ab)^2=1\rangle$, with $n>2$, and take $T=\langle
a\rangle\cong C_n$. Since $T$ is generated by the elements of
order $n$ in $S$, it is a characteristic subgroup of $S$. Now
$A={\rm Aut}\,S$ is isomorphic to the holomorph of $C_n$, a split
extension of $C_n$ by ${\rm Aut}\,C_n=U_n$, the group of units
mod~$(n)$: the normal subgroup $C_n$ is generated by the
automorphism fixing $a$ and sending $b$ to $ab$, while the
complement $U_n$ is the automorphism group of $T$, fixing $b$. If
$n$ is odd, then $Z=1$, so $\overline S=S$ and $\overline T=T$;
since $A/\overline T$ is isomorphic to $U_n$ it is abelian, so
Theorem \theKK\, implies that $S$ cannot be a chirality group. If $n=2m$
is even, then $Z=\langle a^m\rangle\cong C_2$, so $\overline
S\cong D_m$ and $\overline T\cong C_m$; in this case,
$A/\overline T\cong U_n\times C_2$ which is again abelian, so we
have the same result.
\hspace{\fill} $\Box$
\vskip5pt

The above proof fails for $n=2$ since $C_2$ is not characteristic
in $D_2$; indeed, $D_n$ is abelian for $n\leq 2$, and is
therefore a chirality group by Theorem \theJJ.

\vskip8pt
\noindent {\bf Corollary \theteocounter\,\,} {\sl If $q=p^e$ where $p$ is
prime, and $d$
divides $p-1$, then $PGL_d(q)$ and $GL_d(q)$ are not chirality
groups.}
\vskip5pt
\newcounter{OO}
\setcounter{OO}{\theteocounter}
\addtocounter{teocounter}{1}
\noindent{\sl Proof.} Take $S=PGL_d(q)$ and $T=PSL_d(q)$, so
$Z=1$ and $\overline S/\overline T=S/T\cong F_q^*/(F_q^*)^d\cong
C_{(d,q-1)}$. We have $A=P\Gamma L_d(q)$ with $A/T$ abelian
(isomorphic to $C_{(d,q-1)}\times C_e=C_d\times C_e$) if $d$
divides $p-1$, in which case $PGL_d(q)$ cannot be a chirality
group by Theorem \theKK. The same applies to $S=GL_d(q)$, with
$T=SL_d(q)$ and $Z$ the group of scalar matrices, so that
$\overline S=PGL_d(q)$ and $\overline T=PSL_d(q)$.
\hspace{\fill} $\Box$
\vskip5pt

In particular, if $q$ is odd then $PGL_2(q)$ and $GL_2(q)$ are not
chirality groups.

\vskip8pt

\noindent{\sl Proof of Theorem~\theKK.} Suppose that $S$ is a
chirality group, so $\Delta^+$ has a normal subgroup $H$ with
$H/H_{\Delta}\cong H^r/H_{\Delta}\cong H^{\Delta}/H\cong
H^{\Delta}/H^r\cong S$, in the usual notation. Let
$Q=\Delta^+/H_{\Delta}$, so $Q$ has normal subgroups
$N_1=H/H_{\Delta}$ and $N_2=H^r/H_{\Delta}$, both isomorphic to
$S$, generating their direct product $N_1\times N_2\cong S\times
S$ as a normal subgroup $N=H^{\Delta}/H_{\Delta}$ of $Q$. Now let
$C$ be the centraliser $C_Q(N)$, a normal subgroup of $Q$, and
for each subgroup $X\leq Q$ let $\overline X=XC/C\cong X/(X\cap
C)$, the image of $X$ in $Q/C$. We have $N\cap C=Z(N_1)\times
Z(N_2)\cong Z(S)^2$, so $\overline N=\overline N_1\times
\overline N_2\cong (S/Z(S))^2$, with $\overline N$ and each
$\overline N_i$ normal in $\overline Q$. The action of $Q$ by
conjugation on $N$, preserving each $N_i$, induces a faithful
action of $\overline Q$; this gives an embedding $\overline
Q\leq{\rm Aut}\,N_1\times{\rm Aut}\,N_2\cong A^2$, where $A$
denotes ${\rm Aut}\,S$.

In the natural homomorphism $S\to \overline S=S/Z(S)={\rm
Inn}\,S\leq A$, where $S$ acts on itself by conjugation, any
characteristic subgroup $T$ of $S$ maps to a normal subgroup
$\overline T$ of $A$, contained in $\overline S$. Suppose that
$\overline T<\overline S$ and $A/\overline T$ is abelian. If
$\overline T_i$ is the subgroup of $\overline N_i$ corresponding
to $\overline T$, then $\overline T_1\times\overline N_2$ is a
normal subgroup of $\overline Q$, with
$$\overline Q/(\overline T_1\times\overline N_2)\leq A^2/(\overline
T_1\times\overline N_2)\cong A/\overline T_1\times A/\overline
N_2.$$ Now $A/\overline T_1$ is abelian, as is its quotient
$A/\overline N_2$, and hence so is $\overline Q/(\overline
T_1\times\overline N_2)$. The automorphism $\alpha$ of $\Delta^+$
acts on $Q$, preserving $N$ and hence $C$, so it induces an
automorphism of $\overline Q$, inverting its two generators;
since $\overline Q/(\overline T_1\times\overline N_2)$ is
abelian, it has an automorphism inverting its two generators, so
$\alpha$ preserves $\overline T_1\times\overline N_2$. However,
it transposes $\overline N_1$ and $\overline N_2$ and hence
transposes $\overline T_1\times\overline N_2$ and $\overline
N_1\times\overline T_2$, which are distinct since $\overline
T_i<\overline N_i$. This contradiction shows that $S$ cannot be a
chirality group.
\hspace{\fill} $\Box$
\vskip5pt

A similar argument shows that $A_5$ cannot be a chirality group.
For otherwise, putting $S=A_5$ we have $Z(S)=1$ and $A=S_5$, so
$A_5\times A_5\leq\overline Q\leq S_5\times S_5$ with $\alpha$
inducing an automorphism of $\overline Q$ which transposes the
two copies of $A_5$. This implies that $\Delta^+$ has a normal
subgroup, with quotient group $A_5$ or $S_5$, which is not normal
in $\Delta$, contradicting Corollary \theCC.

\vskip8pt

In contrast with the rather simple criteria in Theorem \theJJ\, and
Corollaries \theMM\, and \theNN\, the next result suggests that a complete
characterisation of chirality groups will be a difficult task.

\vskip8pt
\noindent {\bf Theorem \theteocounter\,\,} {\sl The 1-dimensional affine
group $AGL_1(q)$ is
a chirality group if and only if $q$ is an odd power of
2}
\vskip5pt
\newcounter{PP}
\setcounter{PP}{\theteocounter}
\addtocounter{teocounter}{1}
\noindent{\sl Proof.} If $S=AGL_1(q)$ then $Z=Z(S)=1$, and we can
identify $A={\rm Aut}\,S$ with $A\Gamma L_1(q)$. This group
consists of the transformations $t\mapsto at^{\gamma}+b$ of $F_q$,
where $a,b\in F_q, a\neq 0$ and $\gamma\in\Gamma={\rm Gal}\,F_q$;
this Galois group is a cyclic group of order $e$ generated by the
Frobenius automorphism $t\mapsto t^p$ of $F_q$, where $q=p^e$ and
$p$ is prime. The transformations with $\gamma=1$ form the group
$S$, a normal subgroup of index $e$ in $A$ on which $A$ acts by
conjugation.

In the proof we will use the following subgroups of $S$: for each
divisor $d$ of $q-1$ let $T_d$ be the group of transformations
$t\mapsto at+b$ where $a$ is a $d$-th power in $F_q^*$. This is
the unique subgroup of index $d$ in $S$, so it is characteristic
in $S$; $A/T_d$ is a split extension of $S/T_d\cong C_d$ by
$A/S\cong\Gamma\cong C_e$.

If $q$ is odd, we can apply Theorem \theKK\, with $T=T_2$. We have
$\overline S=S$ and $\overline T=T$, so $\overline T<\overline
S$. Since $A/\overline T=A/T\cong C_2\times C_e$, which is
abelian, Theorem \theKK\, implies that $S$ cannot be a chirality group.

Now let $q=2^e$. First we show that if $e$ is odd then $S$ is a
chirality group. Since $AGL_1(2)\cong C_2$, we may assume that
$e>1$. Let $x\in A$ be given by $t\mapsto t^2$, and let $y\in S$
be given by $t\mapsto u(t-v)+v=ut+(1-u)v$, where $u$ generates
$F_q^*$ and $v\in F_q\setminus F_2$. Then $y$ generates the
stabiliser $S_v$ in $S$ of $v$, a cyclic group of order $q-1$
consisting of the transformations $t\mapsto u^i(t-v)+v$, and its
conjugate $y^x=x^{-1}yx:t\mapsto (u\sqrt
t+(1-u)v)^2=u^2t+(1-u)^2v^2$ generates the stabiliser $S_{v^2}$ of
$v^2$. Since $v\neq 0,1$ we have $v\neq v^2$, so $S_v\neq
S_{v^2}$. Since $S$ acts primitively (in fact, doubly
transitively) on $F_q$, the stabilisers of points are maximal
subgroups of $S$, so $S=\langle S_v, S_{v^2}\rangle=\langle y,
y^x\rangle$; since $xS$ generates $A/S$ it follows that
$A=\langle x,y,y^x\rangle=\langle x,y\rangle$. This gives an
epimorphism $\theta:\Delta^+\to A$, so that $A$ is the monodromy
group of a regular orientable hypermap, corresponding to the
normal subgroup $H={\rm ker}\,\theta$ of $\Delta^+$.

The normal subgroup $K=\theta^{-1}(S)$ of $\Delta^+$ satisfies
$\Delta^+/K\cong A/S\cong{\rm Gal}\,F_q\cong C_e$, so
$\Delta^+/K$ is abelian; it follows from Corollary \theCC\, that $K$ is
normal in $\Delta$ and hence $H^{\Delta}\leq K$. If we can show
that $H^{\Delta}=K$ then we have a chirality group
$H^{\Delta}/H=K/H\cong S$, as required, so suppose that
$H^{\Delta}<K$. Since $H^{\Delta}\geq H$ we have
$H^{\Delta}=\theta^{-1}(N)$ for some normal subgroup $N$ of $A$,
properly contained in $S$, and hence contained in a maximal
normal subgroup of $S$. Now the maximal normal subgroups of $S$
are the subgroups $T_p$ defined earlier, one for each prime $p$
dividing $q-1$, so $N\leq T_p$ for some $p$. Since $q$ is even,
$p$ is odd. We have
$$A/T_p=\langle x,y\mid x^e=y^p=1, y^x=y^2\rangle,$$
so in $A/T_p$ we have
$$(y^{-1})^{x^{-1}}=y^{-1/2}\;(=y^{(p-1)/2}).$$
Now the automorphism $\alpha$ of $\Delta^+$ leaves $H^{\Delta}$
and $K$ invariant, so it induces an automorphism of
$A/N\;(\cong\Delta^+/H^{\Delta})$ leaving $S/N\;(\cong
K/H^{\Delta})$ invariant; since $T_p/N$ is the unique normal
subgroup of index $p$ in $S/N$ it is invariant under $\alpha$,
which therefore induces an automorphism of $A/T_p$. This
automorphism inverts $x$ and $y$, so by applying it to the
relation $y^x=y^2$ of $A/T_p$ we get
$$(y^{-1})^{x^{-1}}=y^{-2}$$
in $A/T_p$. Comparing these two equations we get
$y^{-1/2}=y^{-2}$, so $y^3=1$ in $A/T_p$ and hence $p=3$. Thus
$q=2^e\equiv 1$ mod~$(3)$, so $e$ is even, against our
hypothesis. This contradiction shows that $AGL_1(2^e)$ is a
chirality group for all odd $e$.

If $e=2f$ is even then $q\equiv 1$ mod~$(3)$, so $S$ has a
characteristic subgroup $T=T_3$ with
$$A/T=\langle x,y \mid x^e=y^3=1, y^x=y^2\rangle$$
and $S/T=\langle y\rangle\cong C_3$. The subgroup $\langle
x^2\rangle$, isomorphic to $C_f$, is normal (in fact, central) in
$A/T$, with quotient isomorphic to $D_3$; let $U$ denote its
inverse image in $A$, so $T\leq U\leq A$ with $U/T\cong C_f$ and
$A/U\cong D_3$. For any subgroup $X$ of $A$ (or of $A^2$), let
$X^*$ denote its image $XU/U\cong X/(X\cap U)$ in $A/U$ (or in
$A^2/U^2$), so $S^*\cong C_3$ and $A^*\cong D_3$. We have
$T^2<S^2\leq\overline Q\leq A^2$, so $C_3^2\cong
(S^2)^*\leq\overline Q^*\leq(A^2)^*\cong D_3^2$ and hence
$\overline Q^*$ is either $C_3^2$ or $D_3^2$ or an extension of
$C_3^2$ by $C_2$. Now $\alpha$ transposes $H$ and $H^r$, so it
induces an automorphism of $\overline Q^*$, inverting its
generators and transposing the corresponding two normal subgroups
$C_3$. However, these two normal subgroups of $\overline Q^*$
have quotient group $C_3$ or $D_3$ or $D_3\times C_2\cong D_6$,
which are abelian or dihedral, so in each case Corollary
\theBB\,
implies that they are invariant under $\alpha$. This
contradiction shows that $AGL_1(2^e)$ is not a chirality group
when $e$ is even.
\hspace{\fill} $\Box$


\vskip5pt

\vskip10pt

\begin{thebibliography}{99}

\bibitem{BAR} L.~D.~Barron, Fundamental symmetry aspects of molecular
chirality,
in {\sl New Developments in Molecular Chirality} (P.~G.~Mezey ed.), Kluwer
Acad.~Pub.,
Dordrecht, 1991.

\bibitem{Bre} A.~J.~Breda d'Azevedo and G.~A.~Jones, Double coverings and
reflexible abelian hypermaps, {\sl Beitr\"age zur Algebra und Geometrie} 41 (2000), 371--389.

\bibitem{CHIF} A.~Breda d'Azevedo and R.~Nedela,  Chiral hypermaps with few
hyperfaces, Math. Slovaca, {\bf 53}, No. 2, (2003), 107-128.

\bibitem{Bro} C.~Brown (ed.), {\sl Chirality in Drug Design and Synthesis},
Academic Press, London,
1990.

\bibitem{CCNPW} J.~H.~Conway, R.~T.~Curtis, S.~P.~Norton, R.~A.~Parker,
R.~A.~Wilson, {\sl ATLAS of Finite Groups}, Clarendon Press, Oxford,
1985.

\bibitem{Cori} R.~Cori and A.~Machi, Maps, hypermaps and
their automorphisms: a survey I, II, III {\sl Expositiones
Math.} 10 (1992), 403--427, 429--447, 449--467.

\bibitem{CM} H.~S.~M.~Coxeter and W.~O.~J.~Moser, {\sl Generators and
Relations for Discrete
Groups}, 3rd ed., Springer-Verlag, New York, 1972.

\bibitem{CS} D.~Corn and D.~Singerman, Regular hypermaps, {\sl
Europ.~J.~Comb.} 9 (1988),
337-351.

\bibitem{FUJ} S.~Fujita, {\sl Symmetry and Combinatorial Enumeration in
Chemistry}, Springer-Verlag, Berlin, 1991.

\bibitem{GKS} R.~Guralnick, W.~Kantor and J.~Saxl, Probability of
generating a classical
group, {\sl Comm.~Alg.\/} 22 (1994), 1395--1402.

\bibitem{HB} B.~Huppert and N.~Blackburn, {\sl Finite Groups III\/}, Springer,
Berlin, 1982.

\bibitem{JAN} R.~Janoschek (ed.), {\sl Chirality -- From Weak Bosons to
$\alpha$-Helix},
Springer-Verlag, Berlin, 1991.

\bibitem{Jon} G.~A.~Jones, Ree groups and Riemann surfaces, {\sl
J.~Algebra} 165
(1994), 41--62

\bibitem{JS} G.~A.~Jones and S.~A.~Silver, Suzuki groups and surfaces, {\sl
J.~London
Math.~Soc.\/} (2) 48 (1993), 117--125.

\bibitem{KEL} Lord Kelvin, {\sl Baltimore lectures}, C.~J.~Clay, London, 1904.

\bibitem{Mac} A.~M.~Macbeath, On a theorem of Hurwitz, {\sl Proc.~Glasgow
Math.~Assoc.} 5 (1961), 90--96.

\bibitem{MEZ} P.~G.~Mezey (ed.), {\sl New Developments in Molecular
Chirality\/},
Kluwer Acad.~Pub., Dordrecht, 1991.

\bibitem{Re1} R.~Ree, A family of simple groups associated with the simple
Lie algebra of type $(G_2)$, {\sl Bull.~Amer.~Math.~Soc.} 66 (1960),
508--510.

\bibitem{Re2} R.~Ree, A family of simple groups associated with the simple
Lie algebra of type $(G_2)$, {\sl Amer.~J.~Math.} 83 (1961), 432--462.

\bibitem{Sin} D.~Singerman, Symmetries of Riemann surfaces with large
automorphism group, {\sl Math.~Ann.} 210 (1974), 17--32.

\bibitem{Suz}  M.~Suzuki, On a class of doubly transitive groups, {\sl
Ann.~Math.}
(2) 75 (1962), 105--145.

\bibitem{Wa} T.~R.~S.~Walsh, Hypermaps versus bipartite maps, {\sl J.
Combinatorial Theory (B)} 18 (1975), 155--163.

\bibitem{Wie} H.~Wielandt, {\sl Finite Permutation Groups\/}, Academic
Press, New
York, 1964.
\end{thebibliography}
\end{document}